\documentclass[11pt]{article}

\usepackage{amsmath}
\usepackage{amssymb}
\usepackage{amsthm}
\usepackage{dsfont}
\usepackage{amsfonts}
\usepackage{cases}
\usepackage{graphicx}
\usepackage{enumitem}
\usepackage{url}
\theoremstyle{plain}

\theoremstyle{definition}

\newcommand{\R}{\mathbb R}

\newcommand{\E}{\mathrm E}
\newcommand{\F}{{\cal F}}
\newcommand{\Ft}{{\cal F}_t}

\usepackage[longnamesfirst]{natbib}
\RequirePackage[colorlinks,citecolor=blue,urlcolor=blue]{hyperref}

\begin{document}

\title {Last-Hitting Times and Williams' Decomposition \\of the Bessel Process of Dimension 3 \\at its Ultimate Minimum}
\date{}
\author{F. Thomas Bruss and Marc Yor \\Universit\'e Libre de Bruxelles and Universit\'e de Paris VI}
\maketitle

\begin{abstract}\noindent  In this note we shortly recall the importance of last-hitting times in  theory and applications of optimal stopping. As a small contribution to this domain we then propose a concise proof of David Williams' decomposition of the Bessel Process of dimension 3 (BES(3)), starting from $r>0$ at its ultimate minimum. This discussion is strongly motivated by our interest in properties of last hitting times in general, and here specifically, directly linked with the forthcoming reading guide of Nikeghbali and Platen
on this subject.

best choice problems, last arrival problem, compassionate-use clinical trials.


\smallskip \noindent {\bf Running title:} 
Last Hitting Times and Williams' Decomposition
\end{abstract}

\section{Introduction}\label{intro}

 Only in trivial cases, last hitting times are at the same time stopping times because the "current- measurability" requirement is usually not satisfied. Hence it is typically harder to deal with last hitting times than with stopping times for which we know quite an impressive collection of Theorems and tools. As Chung (see citation of Nikeghbali and Platen(2012)) among others conclude, last hitting times must therefore be {\it avoided at all costs.} 

This is one way to see things, but one should admit that, often enough, reality looks somewhat different. Indeed, ironically, many interesting problems in the theory of optimal stopping require us to deal with last hitting times, and not with stopping times. And so, the attitude has changed, and the work of Jeulin (1980) and others had quite an influence on this development. 

In their recent paper, Nikeghbali and Platen cite several interesting examples from the domain of Mathematical Finance. In order to add to the motivation, we would like to slightly broaden the horizon and look at  a few other examples. 

\medskip
 {\bf Examples and Motivation}
 
 \smallskip
Best-choice problems, or secretary problems, are typical representatives of a last-hitting time problem, namely problems of stopping on the last "improvement" of a stochastic process. In some of these problems, the difficulty stemming from the last-hitting time character disappears. To give a very simple example, suppose we observe sequentially variables $X_1, X_2, \cdots$ and would like to maximize, for a given objective function $f,$  the expected total return, that is we seek
$$ \arg \max_\tau f(X_1, X_2, \cdots , X_\tau).$$ Suppose now that the optimal payoff for stopping {\it after} time $t$ does not depend on $\Ft,$ where $(\F_s)$ denotes the natural filtration. Then 
$$\sup_{\tau \ge t} \E\left(f(X_1, X_2, \cdots, X_\tau |\Ft)\right)=\sup_{\tau \ge t}\E(f(X_{t+1}, X_{t+2}, \cdots, X_\tau))$$ so that RHS as well as  $X_1, X_2, \cdots, X_t$ are both $\Ft$-measurable. Hence it suffices to compare at each time $t$
the value $f(X_1, X_2, \cdots, X_t)$ with the RHS supremum  in order to take the optimal decision.

\bigskip
{\bf Harder last-hitting time problems}

\smallskip
In more difficult problems the $\Ft$-independence is  typically no longer satisfied. However,  external information about the underlying process  may help us to change nevertheless the last-hitting time problem into a tractable stopping problem. See for example the {\it last-arrival} problem (Bruss and Yor (2012)) which is a continuous time problem, and where the relevant external information about the underlying process is implied by a related martingale. There may be other examples of such an approach.

There are also certain other problems where the {\it last} hitting time objective is  hiding behind other objectives, as for example the objective to discover the {\it first} time a random subset of a given set becomes complete.
We give one specific example of this in the important field of clinical trials, more precisely, in so-called compassionate-use clinical trials. 

In such trials, a sequence of patients is treated with a drug (sometimes without FDA-approval) which may have serious side effects, the only justification being that it may be the last hope.  (These trials require a special written consent of patients.) Treatments are typically sequential so that the physician or statistician may learn form preceding observations.

 A little reflection shows that these trials  pose a difficult ethical problem. The conscientious physician should try to save all lives which can be saved, and, at the same time avoid all unnecessary sufferings caused by the treatment. Since he or she is no prophet, the goal must be to stop (in a given sequence of patients within a fixed horizon) with maximum probability with the first
patient completing the random subset of successes that is stop  with the last success. Indeed, then {\it all} successes are covered, whereas the remaining patients (de facto not savable by the drug) do not have to suffer unnecessarily.
If the success probability for each patient  is known beforehand , then the optimal strategy follows immediately from the odds-algorithm (see Bruss (2000)). However, the physician may have almost no information about the respective success probabilities, and then the general solution of the optimal stopping problem is an open problem.

\smallskip
These examples second the motivation of the studies of Nikeghbali and Platen, as well as ours in the present note.  Last-hitting times are often important. Admittedly, no model in these papers seems tractable enough to deal for instance with the {\it general} compassionate-use stopping problem described above. However, this indicates that it may be worth trying to view last-hitting times from any possible angle, and this is what we try in this paper. 

Our goal is to add to the understanding by looking at what can be done with enlargements of filtrations (in this case together with Girsanov's Theorem), even though we confine our interest  to special processes.

\bigskip
{\bf The BES(3)-process and Williams' Theorem}

\medskip

\noindent {\bf (1.1)}~~In their survey about last passage times, the Nikeghbali and Platen (2012) illustrate some of their formulae with the following example:

\smallskip
Let $(R_t)_{t\ge 0}$ be a BES(3)-process on $\R_+$ starting from $r>0.$ Denote by $(\Ft)_{t \ge 0}$ its natural filtration, and let $I_t$ denote the current infimum of the process $R$ at time $t$, that is,
\begin{align*} I_t=\inf_{s \le t}R_s.\end{align*}
The following results can be found in Nikeghbali and Platen around Corollary 4.10:

\medskip
(a) $I_{\infty}$ follows the same distribution as the random variable $rU,$ where $U$ is uniform on $[0,1].$

\medskip
(b)  The Az\'ema-supermartingale associated with the random time $g$ at which the process $(R)$ reaches $I_\infty$ is given by \begin{align*}Z_t \equiv P(g>t|\Ft)=\frac{I_t}{R_t}.\end{align*}

\medskip
(c) The Laplace transform of the law of $g$ is
\begin{align*} \E \left(e^{-\lambda g}\right)=
\frac{1}{\sqrt{2 \lambda} r} \left(1-e^{-{\sqrt{2 \lambda} r}}\right).
\end{align*}

\medskip
(d) The density of $g$ denoted by $p(t)$ equals
\begin{align*} p(t)=\frac{1}{\sqrt{2 \pi t}\, r} \left(1-e^{-(r^2/2t)}\right). \end{align*}

\bigskip

Our aim in the remaining part of this note is to show Williams' decomposition of a BES(3)-process at its ultimate minimum, and how this decomposition is closely connected with (a)-(b)-(c)-(d).
\bigskip

\noindent
{\bf (1.2)}~~Recall that if $(B_t)_{t\ge 0}$ is a Brownian motion starting from $0$ and $a$ is a real constant, then the law of the first hitting time of $a$ by $(B_t),$ denoted by $T_a^{(B)}$ is given by
\begin{align} 
P\left(T_a^{(B)}\in dt \right)= \frac{dt}{\sqrt{2 \pi t^3}}\, |\,a\,|\, {\rm exp} \left(-\frac{a^2}{2t}\right).
\end{align}This well-known fact allows us to rewrite the statements (c) and (d) above as 
\begin{align} g \overset{ \cal L}{=} T_{r U}^{(B)}, \end{align}
where  $U$ is independent of $(B)$ and uniform on $[0,1],$ and where  $ \overset{ \cal L}{=}$ denotes identity in law.  This can be verified  using (c) and (d). In fact, (2) may be understood via the classical decomposition of the process $(R)$ before and after time $g,$ due to Williams (1974).

\section{Williams' decomposition of $(R)$, before and after $g$, via progressive enlargement}

\medskip
\noindent {\bf (2.1)}~~Figure 1 describes well the decomposition of a BES(3)-process.
\begin{figure}[h]
\begin{center}
{\includegraphics[width= 0.9\linewidth, height=11cm]{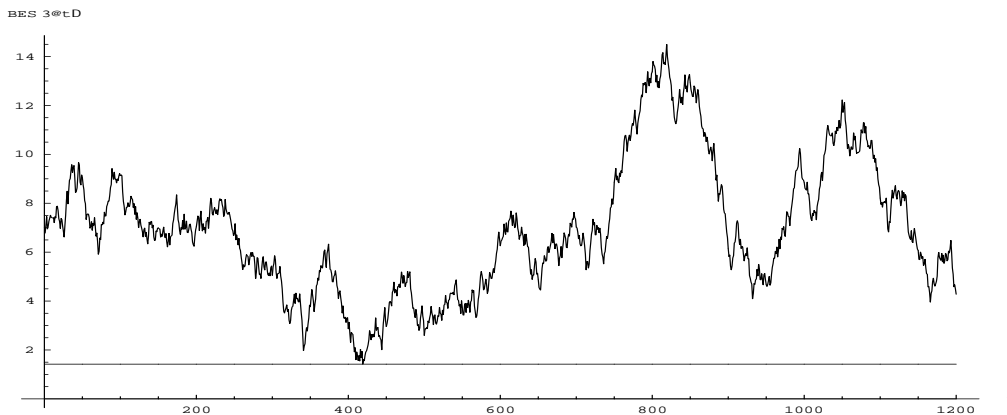}}
 \caption{ This graph presents a simulation of a Bessel process of dimension 3 based on three independent simulations of $U[-1/2,1/2]$-random walks $S^x_k, S^y_k,S^z_k,$  where $k$ runs from $1$ to $1200.$ The starting point is chosen $S^x_0=4, S^y_0=4,S^z_0=2,$ so that the starting level of the simulated process is $B_0=r=6.$ The minimum height is indicated by the supporting horizontal line. In this simulation it equals $1.41963.$ Note that this is $I_{1200}$ and not $I_\infty$, of course. Its level would be uniformly distributed on $[0,B_0]=[0,6]$ if the horizon were infinite. }
\end{center}
\end{figure} 

\noindent Note that this figure is nothing else but a (simulated finite-horizon) version of the Figure 5 in Revuz-Yor (1999) (see Proposition 3.10 and Theorem 3.11 in Ch. 6, Sect. 3) where the BES(3)-process is considered starting from level $c:=r.$

\medskip
\noindent {\bf (2.2)}~~ We now state precisely Williams' Decomposition Theorem before and after time $g.$

\medskip
\noindent {\bf Theorem 1} (Williams (1974))
 
\noindent Consider the following three independent random objects:

\medskip
(i)\,\, ~a Brownian motion $(B'_t)_{t\ge 0}$ with $B'_0=r>0;$

\smallskip
(ii)\, ~a uniform random variable $U$ on $[0,1];$

\smallskip
(iii) ~a BES(3)-process $(\tilde R_t)_{t\ge 0}$ with $\tilde R_0=0;$

\medskip \noindent  Then the process $(R)$ defined by
 \begin{align}
R_t=\begin{cases} B'_t &\mbox{,~if~} t \le g\\  r U+\tilde R_{t-g} & \mbox{,~if~} t \ge g\end{cases}
 \end{align}
with $g=\inf\{u \ge 0: B'_u=rU\}$ is a BES(3)-process  starting from $r>0.$

\bigskip
\noindent
We note that the pre-$g$-Browian-motion 
found in (3) explains the result (2). Indeed, if $B'_t=r-B_t^{(0)}$ then
$$ g=\inf\left\{u\ge 0: B_u^{(0)}=r(1-U)\right\},\eqno(2')$$ which implies (2).

\medskip

 \medskip \noindent
 {\bf 2.3}~We now proceed to the proof of the Theorem via the enlargement formula which describes the additive decomposition of the BES(3)-process $(R_t)$ in the filtration $\left(\F_t^g)\right)$ containing
 the filtration $(\F_t),$ and making $g$ a stopping time.

Firstly, we have $$R_t=r+B_t+\int_0^t \frac{ds}{R_s}, \eqno(4)$$ where $(B_t)$ is a Brownian motion with respect to the filtration $(\F_t).$ 

Secondly, the enlargement formula (see e.g. Jeulin (1980)) yields $$r+B_t=B'_t+\int_0^{g\wedge t} \frac{d<B,Z>_u}{Z_u}+\int_g^t \frac{d<B,1-Z>_u}{1-Z_u} \eqno(5) $$ with $(B'_t)$ being a Brownian motion with respect to $(\F_t^{g}).$

Thirdly, we deduce from (b) the two identities $$\frac{d<B,Z>_u}{Z_u}=-\frac{du}{R_u}, \mbox{~for~} u\le g \eqno(6)$$ and $$\frac{d<B,1-Z>_u}{1-Z_u}=\frac{I_\infty du}{R_u(R_u-I_\infty)} \mbox{,~for~} u>g .\eqno(7)$$ 
These two identities imply (using (4) and (6)) , and also
$$\frac{1}{R_u}+\frac{I_\infty}{R_u(R_u-I_\infty)} =\frac{1}{R_u-I_\infty}$$
the form of the pre-$g$-process  and the form of the post-$g$-process. 

\smallskip
Finally, for the proof of (3) to be complete, it remains to prove that the process $(B')$
is independent of the random variable 
$I_\infty \overset{\cal L}=\,rU,$ or more precisely, that, given 
$I_\infty = a$, the pre-$g$-process 
is just the process $(B'_u)_{u\le T'_a}$ with obvious notation. This is asserted in the following proposition:

\bigskip
{\bf Proposition}: Let $(\Phi_u)_{u\ge 0}$ be a non-negative predictable process on path-space. Further, let $P_r$ denote the law of the process $(R)$ starting from $r$ and let $P'_r$ denote the law of the Brownian motion $(B')$ starting from $r.$ Then, for $a<r,$
\begin{align*}\E_r\left[\Phi_g | I_\infty=a\right]&=\E_r\left[\Phi_{T_a}|T_a <\infty\right]~~~~~~~~(8.1)\\&\equiv \E_r \left[\Phi_{T_a}|I_\infty<a\right]~~~~~~~~~(8.2)\\&=\E'_r\left[\Phi(B'_u; u \le T'_a )\right]~~~~~~(8.3)\end{align*}

\medskip
{\bf Proof}: The equality between the RHS of (8.1) and (8.3) follows, as we will show, from Doob's absolute continuity relationship, namely$$P_r/ \F_t=\left(\frac{X_{t \wedge T_0}}{r}\right)P'_r / \F_t,$$ on the canonical path-space $C([0, \infty],\R),$ where $(X_t)$ denotes the coordinate process on path-space. Indeed, this equality may be extended when replacing the time $t$ by a stopping time. Restricting $\F_{T_a}$ on the set $\{T_a<\infty\}$ we get then in particular $$P_r/(\F_{T_a}\cap \{T_a<\infty\})=\left(\frac{a}{r}\right) P'_r/\F_{T_a}, ~0<a<r,$$ which yields the desired result.

\noindent  identity (8.2) is obvious, since the equality $\{T_a<\infty\}=\{I_\infty<a\}$ holds $P_r$-almost surely.

\medskip
The proof of the equality (8.1) is slightly more subtle.
We start with the identity
$$\E\left[{\bf 1}_{\{g \le t\}}\varphi(I_\infty)\right]= E\left[(1-Z_t)\,\int_0^t\varphi(I_s)d(1-Z_s)\right]\eqno(9)$$
which holds for any Borel-measurable function $\varphi:[0, \infty[ \to \R_+.$ To see this, note that
$$\E\left[{\bf 1}_{\{g\le t\}}\varphi(I_t)\right]=\E\left[(1-Z_t)\varphi(I_t)\right].$$ Assuming $\varphi \in {\cal C}^1$ the latter becomes by partial integration $$\E\left[\int_0^t\varphi'(I_s)dI_s(1-Z_s)\right]+\E\left[\int_0^t\varphi(I_s)d(1-Z_s)\right].$$ We note that the expectation involving $\varphi'$ vanishes, since $1-Z_s$ vanishes $dI_s$ almost everywhere. Thus a monotone class argument implies that (9) holds for every non-negative Borel-measurable function $\varphi.$

\smallskip
Next, from the additive decomposition of $(1-Z_s)$, we obtain $$\E\left[{\bf 1}_{\{g\le t\}}\varphi(I_\infty)\right]=\E\left[\int_0^t\varphi(I_s)\left(-\frac{dI_s}{I_s}\right)\right]=\E\left[\int_{I_t}^r\varphi(a) \frac{da}{a}\right].$$  
Since $\{I_t\le a\}=\{t \ge T_a\}$ the latter can also be written as 
$$\E\left(\int_0^r \varphi(a)\frac{da}{a}{\bf 1}_{\{T_a\le t\}}\right),$$ so that from (9) 
$$\E\left[{\bf 1}_{\{g\le t\}}\varphi(I_\infty)\right]=\E\left(\int_0^r \varphi(a)\frac{da}{a}{\bf 1}_{\{ T_a \le t\}}\right) \eqno(10).$$ Now note that the identity (10) still holds if we replace $t$ by a generic stopping time. Applying again the monotone class theorem gives us then $$\E\left[\Phi_g\varphi(I_\infty)\right]=\E\left[ \int_0^r \varphi(a) \frac{da}{a}\Phi_{T_a}{\bf 1}_{\{T_a < \infty\}}\right]. \eqno (11)$$
Finally, using $I_\infty \overset{{\cal L}}= rU$ with $U$ being uniform on $[0,1]$ under $P_r$ (see (2)), we see that identity (11) implies identity (8.1).

\qed
\section{Concluding remarks}
The statement of the theorem invites for a proof chosing between, on the one hand, initial enlargement with $I_\infty,$ and, on the other hand,  progressive enlargement with $g.$ However we have not exactly proceeded like this; the Proposition plays the role of the initial enlargement method and relies on a classical Girsanov relationship between $P_r$ and $P'_r$. 

In conclusion we find it interesting to present the above as an example of the potential of enlargement techniques,
and here specifically, of a melange of enlargement techniques and Girsanov's theorem. Having said so, we know that this approach can in priciple be done for higher dimensions; however, the corresponding Theorem 1 would look more complicated.

\bigskip

\smallskip

{\bf Acknowledgement}

\smallskip
The authors are grateful to Monique Jeanblanc for providing them with a preprint of Nikeghbali and Platen (to appear).

\bigskip

\smallskip
{\bf References}:

\smallskip
F.T. Bruss: {\it Sum the odds to one and stop}, Ann. Probab., Vol. 28, No 3, 1384-1391, (2000)

\smallskip
F.T. Bruss and M. Yor: {\it Stochastic processes with proportional increments and the last-arrival problem},
Stoch. Proc. and Th. Applic.,  Vol. 122 (9),  3239-3261, (2012)

\smallskip
Th. Jeulin: {\it Semi-martingales et grossissement d'une filtration}, LNM Springer 833, (1980).

\smallskip
A. Nikeghbali and E. Platen: {\it A reading guide for last passage times with financial applications in view}. To appear in Finance and Stochastics (2012)

\smallskip
D. Revuz and M. Yor: {\it Continuous martingales and Brownian motion.}3rd Edition, Springer (1999).

\smallskip
D. Williams: {\it Path decomposition and continuity of local times for one-dimensional diffusions, I}, Proceed. of London Math. Soc. (3), 28, 738-768, (1974).
\bigskip

\bigskip
\noindent Authors' adresses:

\bigskip

\noindent
Universit\'e Libre de Bruxelles\\Facult\'e des sciences\\D\'epartement de Math\'ematique, Campus Plaine CP 210\\ B-1050 Brussels, Belgium.
 
 \bigskip
 \noindent
 Universit\'e Pierre et Marie Curie\\ Laboratoire des Probabilit\'es\\4, place Jusieu, Tour 56, \\F-75252 Paris Cedex 05, France

\end{document}